\documentclass[11pt]{article}
\usepackage{amsmath}
\usepackage{amssymb}
\usepackage{amsthm}
\usepackage[usenames]{color}
\usepackage{amscd}
\usepackage{dsfont}
\usepackage{indentfirst}

\usepackage[colorlinks=true,linkcolor=blue,filecolor=red,
citecolor=webgreen]{hyperref}
\definecolor{webgreen}{rgb}{0,.5,0}
\definecolor{amber}{rgb}{1.0, 0.49, 0.0}
\numberwithin{equation}{section}

\def\N{{\mathds{N}}}

\def\1{{\bf 1}}

\newtheorem{theorem}{Theorem}[section]

\newtheorem{corollary}[theorem]{Corollary}

\begin{document}

\title{{\bf Menon-type identities concerning subsets of the set $\{1,2,\ldots,n\}$}}
\author{L\'aszl\'o T\'oth\thanks{The research was financed by NKFIH in Hungary, within the framework of 
the 2020-4.1.1-TKP2020 3rd thematic programme of the University of P\'ecs.}
\\ Department of Mathematics \\
University of P\'ecs \\
Ifj\'us\'ag \'utja 6, 7624 P\'ecs, Hungary \\
E-mail: {\tt ltoth@gamma.ttk.pte.hu} }
\date{}
\maketitle

\centerline{Mathematica Pannonica New Series, vol. 28 (2022), 65--68}

\begin{abstract} We prove certain Menon-type identities associated with the subsets of the set $\{1,2,\ldots,n\}$ 
and related to the functions $f$, $f_k$, $\Phi$ and $\Phi_k$, defined and investigated by Nathanson
\cite{Nat2007}.

\end{abstract}



\section{Introduction}

Menon's identity states that for every $n\in \N:=\{1,2,\ldots \}$,
\begin{equation} \label{Menon_id}
M(n):= \sum_{\substack{a \text{ (mod $n$)} \\ (a,n)=1}} (a-1,n) = \varphi(n) \tau(n),
\end{equation}
where $a$ runs through a reduced residue system (mod $n$), $(k,n)$ stands for the greatest common divisor (gcd) of $k$ and $n$,
$\varphi(n)$ is Euler's totient function and $\tau(n)=\sum_{d\mid n} 1$ is the divisor function. Identity \eqref{Menon_id}
is due to P. K. Menon \cite{Men1965}, and it has been generalized in various directions by several authors, also in recent papers. 
See, e.g., \cite{CaiMig,CQ2021,Hau2005,HauTot2020,JW2020,Tot2011,Tot2018,Tot2020,WZJ2019B} and their references.
Also see the quite recent survey by the author \cite{Tot2021}.

For a nonempty subset $A$ of $\{1,2,\ldots, n\}$ let $(A)$
denote the gcd of the elements of $A$. Then $A$ is said to be
relatively prime if $(A)=1$, i.e., the elements of $A$ are
relatively prime. Let $f(n)$ denote the number of relatively prime
subsets of $\{1,2,\ldots, n\}$. Here $f(1)=1$, $f(2)=2$, $f(3)=5$, $f(4)=11$, $f(5)=26$, $f(6)=53$; this is sequence A085945 in \cite{Slo}. For every $n\in \N$ one has

\begin{equation} \label{f_relprime}
f(n)=\sum_{d=1}^n \mu(d)\left(2^{\lfloor n/d\rfloor}-1\right),
\end{equation}
where $\mu$ is the M\"obius function and  $\lfloor x \rfloor$ is the floor function of $x$. 

A similar formula is valid for the number $f_k(n)$ of relatively prime
$k$-subsets (subsets with $k$ elements) of $\{1,2,\ldots, n\}$. Namely, for every $n,k\in \N$ ($k\le n$),
\begin{equation} \label{f_relprime_k}
f_k(n)=\sum_{d=1}^n \mu(d) \binom{\lfloor n/d\rfloor}{k}.
\end{equation}

Note that for $k=1$ one has, by a well-known identity (see, e.g. \cite[Eq.\ (2.17)]{Ten2015}),
\begin{equation} \label{f_1}
f_1(n)=\sum_{d=1}^n \mu(d) \lfloor n/d\rfloor = 1  \quad (n\in \N).
\end{equation}

If $k=2$, then $f_2(n)$ is sequence A015614 in \cite{Slo}, namely
\begin{equation*} 
f_2(n)=\sum_{j=2}^n  \varphi(j) \quad (n\in \N),
\end{equation*}
where $f_2(1)=0$, $f_2(2)=1$, $f_2(3)=3$, $f_2(4)=5$, $f_2(5)=9$, $f_2(6)=11$. Also, $f_n(n)=1$ ($n\in \N$).

Furthermore, consider the Euler-type functions $\Phi(n)$ (sequence A027375 in \cite{Slo}) and
$\Phi_k(n)$, representing the number of nonempty subsets $A$ of
$\{1,2,\ldots, n\}$ and $k$-subsets $A$ of $\{1,2,\ldots, n\}$,
respectively, such that $(A)$ and $n$ are relatively prime. Observe that
$\Phi_1(n)=\varphi(n)$ is Euler's function and $\Phi_n(n)=1$ ($n\in \N$). One has $\Phi(1)=1$, 
\begin{equation*} 
\Phi(n)=\sum_{d\mid n} \mu(d) 2^{n/d} \quad (n\in \N, n>1),
\end{equation*}
and for every $n,k\in \N$ ($k\le n$),
\begin{equation*} 
\Phi_k(n)=\sum_{d\mid n} \mu(d) \binom{n/d}{k}.
\end{equation*}

The functions $f$, $f_k$, $\Phi$ and $\Phi_k$ have been defined and investigated by Nathanson
\cite{Nat2007}. Also see the author \cite{Tot2010Integers} and its references.

In this note we present certain Menon-type identities associated with the subsets of the set $\{1,2,\ldots,n\}$,
not investigated in the literature, and related to the above functions.

\section{Results}
We define the sum $\overline{M}(n)$ by
\begin{equation*} 
\overline{M}(n):= \sum_{\substack{\emptyset \ne A \subseteq \{1,2,\ldots,n\} \\ ((A),n)=1}} ((A)-1,n),
\end{equation*}
taken over all nonempty subsets of $\{1,2\ldots,n\}$ such that $(A)$ and $n$ are relatively prime,  
where $((A)-1,n)$ denotes the gcd of $(A)-1$ and $n$. 
Also, for $1\le k\le n$ let
\begin{equation*} 
\overline{M}_k(n):= \sum_{\substack{A \subseteq \{1,2,\ldots,n\} \\ \# A =k \\ ((A),n)=1}} ((A)-1,n),
\end{equation*}
the sum being over the $k$-subsets of $\{1,2,\ldots,n\}$ such that $(A)$ and $n$ are relatively prime.
Observe that the sums $\overline{M}(n)$ and $\overline{M}_k(n)$ have $\Phi(n)$, respectively $\Phi_k(n)$ terms.

If $k=1$, then $\overline{M}_1(n)=M(n)=\varphi(n)\tau(n)$, according to Menon's identity \eqref{Menon_id}.
If $k=n$, then $\overline{M}_n(n)=n$ ($n\in \N$).

We show that for every $n$ and $k$, the values $\overline{M}(n)$ and $\overline{M}_k(n)$ can be expressed as 
linear combinations of the values $f(j)$ ($1\le j \le n$) and  $f_k(j)$ ($1\le j \le n$), respectively. 
More exactly we have the following results.

\begin{theorem} For every $n,k\in \N$,
\begin{equation} \label{overline_M}
\overline{M}(n)=  \sum_{d\mid n} \varphi(d) \sum_{\substack{\delta \mid n\\ (\delta,d)=1}} \mu(\delta) \sum_{\substack{j=1\\ \delta j \equiv 1 
\textup{ (mod $d$)}}}^{n/\delta} 
f\left(\left\lfloor \frac{n}{j \delta}\right\rfloor\right),
\end{equation}
\begin{equation} \label{overline_M_k}
\overline{M}_k(n)=  \sum_{d\mid n} \varphi(d) \sum_{\substack{\delta \mid n\\ (\delta,d)=1}} \mu(\delta) \sum_{\substack{j=1\\ \delta j \equiv 1 
\textup{ (mod $d$)}}}^{n/\delta} 
f_k\left(\left\lfloor \frac{n}{j \delta}\right\rfloor\right),
\end{equation}
where the functions $f$ and $f_k$ are given by \eqref{f_relprime} and \eqref{f_relprime_k}, respectively.
\end{theorem}

Note that if $k=1$, then $f_1(n)=1$ ($n\in \N$) by \eqref{f_1}, and \eqref{overline_M_k} quickly leads to Menon's identity \eqref{Menon_id}.

\begin{corollary} For every prime power $p^t$ \textup{($t \ge 1$)}, 
\begin{equation} \label{M_p_t}
\overline{M}(p^t)=  \sum_{j=1}^{p^t} f\left(\left\lfloor \frac{p^t}{j}\right\rfloor\right)
- \sum_{j=1}^{p^{t-1}} f\left(\left\lfloor \frac{p^{t-1}}{j}\right\rfloor\right)
+ (p-1) \sum_{s=1}^t p^{s-1} \sum_{m=1}^{p^{t-s}} f\left(\left\lfloor \frac{p^t}{1+(m-1)p^s} \right\rfloor\right),
\end{equation}
\begin{equation} \label{M_k_p_t}
\overline{M}_k(p^t)=  \sum_{j=1}^{p^t} f_k\left(\left\lfloor \frac{p^t}{j}\right\rfloor\right)
- \sum_{j=1}^{p^{t-1}} f_k\left(\left\lfloor \frac{p^{t-1}}{j}\right\rfloor\right)
+ (p-1) \sum_{s=1}^t p^{s-1} \sum_{m=1}^{p^{t-s}} f_k\left(\left\lfloor \frac{p^t}{1+(m-1)p^s} \right\rfloor\right).
\end{equation}
\end{corollary}

\begin{corollary} For every prime $p$, 
\begin{equation} \label{M_p}
\overline{M}(p)=  pf(p) -1+ \sum_{j=2}^p f\left(\left\lfloor \frac{p}{j}\right\rfloor\right),
\end{equation}
\begin{equation} \label{M_k_p}
\overline{M}_k(p)=  pf_k(p) - f_k(1)+ \sum_{j=2}^p f_k\left(\left\lfloor \frac{p}{j}\right\rfloor\right).
\end{equation}
\end{corollary}

In particular, $\overline{M}(1)=1$, $\overline{M}(2)=4$, $\overline{M}(3)=16$, $\overline{M}(4)=46$, 
$\overline{M}(5)=134$, $\overline{M}(6)=320$. Also, $\overline{M}_2(1)=0$,
$\overline{M}_2(2)=2$, $\overline{M}_2(3)=9$, $\overline{M}_2(4)=20$, $\overline{M}_2(5)=46$, $\overline{M}_2(6)=66$.

\section{Proofs}

By using the Gauss formula $n=\sum_{d\mid n} \varphi(n)$ ($n\in \N$) and that $G(n):= \sum_{\delta \mid n} \mu(\delta)=0$ for $n>1$ and $G(1)=1$, 
\begin{equation*}
\overline{M}(n) =  \sum_{\substack{\emptyset \ne A \subseteq \{1,2,\ldots,n\} \\ ((A),n)=1}} \sum_{d\mid ((A)-1,n)} \varphi(d) 
=  \sum_{d\mid n} \varphi(d) \sum_{\substack{\emptyset \ne A \subseteq \{1,2,\ldots,n\} \\ ((A),n)=1\\ d\mid (A)-1}} 1
\end{equation*}
\begin{equation*}
=  \sum_{d\mid n} \varphi(d) \sum_{\substack{\emptyset \ne A \subseteq \{1,2,\ldots,n\} \\ (A)\equiv 1 \text{ (mod $d$)}}} \sum_{\delta \mid ((A),n)} \mu(\delta) 
\end{equation*}
\begin{equation} \label{S} 
= \sum_{d\mid n} \varphi(d) \sum_{\substack{\delta \mid n\\ (\delta,d)=1}} \mu(\delta)  \sum_{\substack{\emptyset \ne A \subseteq \{1,2,\ldots,n\} \\ (A)\equiv 1 \text{ (mod $d$)} \\ \delta \mid (A)}} 1, 
\end{equation}
where the condition $(\delta,d)=1$ comes from $(A)\equiv 1$ (mod $d$) and $\delta \mid (A)$.
Also, $\delta \mid (A)$ is equivalent to $A=\delta B:=\{\delta b: b\in B\}$, and we conclude that the last 
sum $S$ in \eqref{S} is 
\begin{align*}
S:= \sum_{\substack{\emptyset \ne A \subseteq \{1,2,\ldots,n\} \\ (A)\equiv 1 \text{ (mod $d$)} \\ \delta \mid (A)}} 1 = \sum_{\substack{\emptyset \ne \delta B \subseteq \{1,2,\ldots,n\} \\ (\delta B)\equiv 1 \text{ (mod $d$)} }} 1 = \sum_{\substack{\emptyset \ne B \subseteq \{1,2,\ldots, n/\delta\} \\ \delta (B)\equiv 1 \text{ (mod $d$)} }} 1.
\end{align*}

Now by grouping the terms of the latter sum according to the values $(B)=j$, where $j=1,2,\ldots,n/\delta$, and denoting $B=jC$ with $(C)=1$ we have
\begin{equation*}
S=  \sum_{\substack{j=1\\ \delta j\equiv 1 \text{ (mod $d$)} }}^{n/\delta}  \sum_{\substack{\emptyset \ne B \subseteq \{1,2,\ldots, n/\delta\} \\ (B)=j}} 1 
= \sum_{\substack{j=1\\ \delta j\equiv 1 \text{ (mod $d$)} }}^{n/\delta}  \sum_{\substack{\emptyset \ne C \subseteq \{1,2,\ldots, \lfloor n/(j\delta) \rfloor \} \\ (C)=1}} 1  \end{equation*}
\begin{equation} \label{S_result}
=  \sum_{\substack{j=1\\ \delta j\equiv 1 \text{ (mod $d$)} }}^{n/\delta} f(\lfloor n/(j\delta) \rfloor), 
\end{equation}
by the definition of the function $f$. Inserting \eqref{S_result} into \eqref{S} the proof of identity \eqref{overline_M} is complete.

The proof of identity \eqref{overline_M_k} is similar.

If $n=p^t$ ($t\ge 1$) is a prime power, then the only nonzero terms in \eqref{overline_M} and \eqref{overline_M_k} are those for $(d,\delta)=(1,1),(1,p),(p,1),(p^2,1),\ldots, (p^t,1)$. This gives \eqref{M_p_t} and \eqref{M_k_p_t}. 

Finally, \eqref{M_p} and \eqref{M_k_p} are obtained from \eqref{M_p_t}, respectively \eqref{M_k_p_t}, in the case $t=1$.

\end{document}